\documentstyle[12pt]{article}

\newcommand{\be}{\begin{eqnarray}}
\newcommand{\ee}{\end{eqnarray}}
\newcommand{\fr}{\frac}
\newcommand{\pa}{\partial}
\newcommand{\intl}{\int\limits}

\def\const {{\rm const}}
\def\epsilon {\varepsilon}
\def\mes {{\rm mes}}
\def\id {{\rm Id}}
\def\b {{\bf B}}

\def\phi {\varphi}

\begin{document}

\begin{center}
{\bf\LARGE On the $L^2$-instability and $L^2$-controllability of steady flows of an
ideal
incompressible fluid}
\end{center}

\bigskip

\centerline{A. Shnirelman}
\medskip
\centerline{School of Mathematical Sciences}
\centerline{Tel Aviv University}
\centerline{69978 Ramat Aviv, Israel}
\centerline{shnirelm@math.tau.ac.il}

\bigskip
\bigskip

{\bf 1. }In this work we are studying the flows
of an ideal incompressible fluid in a
bounded 2-d domain $M\subset{\bf R}^2$,
described by the Euler equations

\be 
&&\frac{\pa u}{\pa
t}+(u,\nabla)u+\nabla p=0;\\
&&\nabla\cdot u=0.  
\ee

Here $u=u(x,t), \ x\in M, \ t\in[0,T]$, and $u|_{\pa M}$
is tangent to $\pa M$.

It has been known for a long time, that if the initial velocity field
$u(x,0)$
is smooth, then there exists unique smooth solution $u(x,t)$ of the Euler
equations, which is defined for all $t\in{\bf R}$; see[M-P]. The next
natural
question is, what may be the behavior of this solution, as $t\to\infty$.
This is a problem of indefinite complexity. A restricted problem is the
following: suppose that the initial flow field $u_0(x)$ is in close to a
steady solution $u_0(x)$. What may happen with this flow for big $t$?
Does it stay always close to $u_0$, or it can escape far away?  Which
flows are available, if we start from different initial velocities, close
to $u_0$? 

These are problems of a global, nonlinear perturbation theory of
steady solutions of the Euler equations. The first idea is to develop
a linear stability theory. The spectrum of a linearized operator is
always symmetric w.r.t. both the real and the imaginary axes, for the
system is Hamiltonian. Therefore we can never prove an asymptotic
stability by the linear method; at best we can prove the absence of a
linear instability, which, in its turn, may be a tricky business.

The true, nonlinear stability of some classes of steady flows was
first proven by V. Arnold (see [A1], [A2], [AK]). He considered a very
strong restriction on perturbation: the perturbation of the vorticity,
$\omega(0)-\omega_0=\nabla\times u(\cdot,0)-\nabla\times u_0(\cdot)$,
should be small in $L^2(M)$. There are three classes of steady flows
which are stable under perturbations small in this sense.  The first
class contains only one flow with constant vorticity; its stability is
obvious. The second and the third classes consist of steady flows,
corresponding to a strong local maximum, resp. minimum, of the kinetic
energy on the leaf of equivortical fields in the space of all smooth
velocity fields in $M$ (see [AK]). 

But this theory breaks down, if we drop the condition on the vorticity
perturbation, and consider all (smooth) velocity fields $u(x,0)$,
which are close to $u_0(x)$ in $L^2$, without any conditions on the
derivatives. Note that this class of perturbations is no less
physically significant than the previous one, because it describes
perturbations with small energy. In this case, there is apparently no
obstacle preventing the flow from going far away from $u_0$.  So, it
is likely that the flow is unstable. (but this does not prove
instability, for there may be some other reasons for stability, like,
for example in the KAM theory.  This is analogous to the situation in
the 3-d Euler equations.  In the 2-d case, the vorticity is
transported by the flow; this means that there exist infinity of
integrals of motion, namely the moments of vorticity.  These integrals
prevent the solution from forming singularities. In the 3-d case, the
vorticity field is transformed by the flow as a frozen-in vector
field, and we can't extract additional integrals of motion from the
vorticity field. This means that we don't know any obstacle to the
formation of a finite-time singularity from a smooth initial flow. But
we have no examples yet of such singularities.  May be, this work may
give some hints.)

{\bf 2. }In this work we consider a weaker stability problem. Consider
the Euler
equations with a nonzero right hand side (i.e. external force): 

\be
&&\frac{\pa u}{\pa t}+(u,\nabla)u+\nabla p=f;\\
&&\nabla\cdot u=0.
\ee

Here $f=f(x,t)$ is a smooth in $x$ vector field, such that $\nabla\cdot
f=0$, and $f(x,t)|_{\pa M}$ is parallel to $\pa M$. Consider the behavior
of $u(x,t)$, if $f$ is small in the following sense:  $\int_0^T\parallel
f(\cdot,t)\parallel_{L^2} dt$ is small, where $[0,T]$ is the time
interval (assumed to be long), where the flow is considered. For example,
if $f$ has the form $f(x,t)=F(x)\delta(t)$, we return to the initial
stability
problem. 

\medskip

{\bf Definition 1.}{\sl Suppose that $u(x_1,x_2)$ and $v(x_1,x_2)$ are
two steady flows. We say that the force $f$ transfers the flow $u$
into the flow $v$ during the time interval $[0,T]$, if the following
is true: if $w(x,t)$ is the solution of the nonhomogeneous Euler
equations (3), (4) with the initial condition $w(x,0)=u(x)$, then
$w(x,T)=v(x)$.}

\medskip

Consider the simplest basic steady flow, namely a parallel flow. Let
$M$ be a strip $0\le x_2\le 1$ in the $(x_1,x_2)$-plane. We restrict
ourselves to the flows having period $L$ along the $x_2$-axis; this
period is the same for all flows that are considered below. Suppose
that the velocity field $u_0(x)$ has the form $(U(x_2), 0)$, where $U$
is a given smooth function (the velocity profile). The original
problem was, for which profiles $U$ the flow $u_0$ is stable. Our main
result is the following

\medskip

{\bf Theorem 1. } {\sl For every nontrivial (i.e. different from
constant) velocity profile $U$ the flow $u_0$ is $L^2$-unstable. This
means that for every function $U(x_2)\ne\const$ there exists
$C>0$, such that for every $\epsilon>0$ the following is true. There
exist $T>0$ and a smooth force $f(x,t)$, defined in $M\times[0,T]$,
such that $\int_0^T\parallel f(\cdot,t)\parallel_{L^2}dt<\epsilon$,
and $f$ transfers the flow $u_0$ during the time interval $[0,T]$ into
a steady flow $u_1$, such that $\parallel u_0-u_1\parallel_{L^2}>C$.}

\medskip

So, the flow may be considerably changed by arbitrarily small force,
provided the time interval is sufficiently long.

Note that if the force $f$ satisfies a stronger condition
$\int_0^T\parallel\omega(\cdot,t)\parallel_{L^2}dt<\epsilon$, where
$\omega=\nabla\times f$ is the vorticity, then, for every Arnold
stable flow $u_0$, the resulting perturbation at time $t$ will be
small, too. 

This theorem is implied by a much stronger assertion.

\medskip

{\bf Theorem 2.}{\sl Suppose that $U(x_2)$ and $V(x_2)$ are two
velocity profiles, such that $\int_0^1 U(x_2)dx_2=\int_0^1
V(x_2)dx_2$, and $\int_0^1\fr{1}{2}|U(x_2)|^2
dx_2=\int_0^1\fr{1}{2}|V(x_2)|^2 dx_2$; let $u_0(x_1,x_2)=(U(x_2),0),
\ v_0(x_1,x_2)=(V(x_2),0)$ be corresponding steady parallel flows
(having equal momenta and energies). Then for every $\epsilon>0$ there
exist $T>0$ and a smooth force $f(x,t)$, such that $\int_0^T\parallel
f(\cdot,t)\parallel_{L^2}<\epsilon$, and $f$ transfers $u$ into $v$
during the time interval $[0,T$]. }

\medskip

This means that the flow of an ideal incompressible fluid is perfectly
controllable by arbitrarily small force.

{\bf 3. } Theorems 1, 2 are proven by an explicit construction of the
flow.

Note first, that if $U_1, U_2, \cdots, U_N$ are velocity profiles, and
Theorem 2 is true for every pair $(U_i, U_{i+1})$ of velocity
profiles, then we can pass from $U_1$ to $U_N$, simply concatenating the
flows connecting $U_i$ and $U_{i+1}$; thus Theorem 2 is true for the
pair $(U_1, U_N)$. Therefore it is enough to construct the sequence of
steady flows with profiles $U_1, \cdots, U_N$, and the intermediate
nonsteady flows connecting every two successive steady ones.

Note also, that it is enough to construct a sequence of
piecewise-smooth flows, for it is not difficult to smoothen them, so
that the necessary force will have arbitrarily small norm in $L^1(0,T;
L^2(M))$.

As a first step, we change the flow with the profile $U=U_1$ by a
piecewise-constant profile $U_2$ with sufficiently small steps; this
may be done by a force with arbitrarily small norm.

Thus, $U_2(x_2)$ is a step function, $U_2(x_2)=U_2^{(k)}$ for
$x_2^{(k-1)}<x_2<x_2^{(k)}, \ k=1, \cdots, K$. Every next profile
$U_i$
is also a step-wise function. We are free to subdivide the steps and
change a little the values of velocity, if these changes are small
enough.

Every flow $u_k$ is obtained by the previous one $u_{k-1}$ by one of
two operations, described in the following theorems.

\medskip

{\bf Theorem 3. }{\sl Let $U(x_2)$ be a step function,
$U(x_2)=U^{(k)}$
for $x_2^{(k-1)}<x_2<x_2^{(k)}$; let $V(x_2)$ be another step
function, obtained by transposition of two adjacent segments
$[x_2^{(k-1)},x_2^{(k)}]$ and $[x_2^{(k)},x_2^{(k+1)}]$. Let
$u(x_1,x_2), \ v(x_1,x_2)$ be parallel flows with velocity profile
$U(x_2), V(x_2)$. Then for every $\epsilon>0$ there exist $T>0$ and a
piecewise-smooth force $f(x,t)$, such that $\int_0^T\parallel
f(\cdot,t)\parallel_{L^2}<\epsilon$, and the force $f$ transfers the
flow $u$ into the flow $v$ during the time interval $[0,T]$.}

\medskip

To formulate the next theorem, remind the law of an elastic collision
of two bodies. Suppose that two point masses $m_1$ and $m_2$, having
velocities $u_1$ and $u_2$, collide elastically. Then their
velocities after collision will be $v_1=2u_0-u_1, \ v_2=2u_0-u_2$,
where $u_0=(m_1 u_1+m_2 u_2)/(m_1+m_2)$ is the velocity of the center
of masses. The transformation $(u_1, u_2)\to(v_1, v_2)$ is called a
{\it transformation of elastic collision}.

\medskip

{\bf Theorem 4. }{\sl Suppose that the profile $U(x_2)$ is like in
Theorem 3, and the profile $V(x_2)$ is equal to $U(x_2)$ outside the
segment $x_2^{(k-1)}<x_2<x_2^{(k+1)}$; on the last segment,
$V(x_2)=v^{(k)}$, if $x_2^{(k-1)}<x_2<x_2^{(k)}$, and
$V(x_2)=v^{(k+1)}$, if $x_2^{(k})<x_2<x_2^{(k+1)}$, where $(v^{(k)},
v^{(k+1)})$ is obtained from $(u^{(k)}, u^{(k+1))}$ by the
transformation of elastic collision, the lengths
$x_2^{(k)}-x_2^{(k-1)}, \ x_2^{(k+1)}-x_2^{(k)}$ playing the role of
masses $m_1, m_1$. Let $u(x_1,x_2), v(x_1,x_2)$ be parallel
flows with profiles $U(x_2), V(x_2)$. Then for every $\epsilon>0$
there exist $T>0$ and a force $f(x,t)$, such that $\int_0^T\parallel
f(\cdot,t)\parallel_{L^2}<\epsilon$, and the force $f$ transfers the
flow $u$ into flow $v$.}

\medskip

Suppose now, that $U(x_2)$ and $V(x_2)$ are two velocity profiles,
having equal momenta and energies. Then it is not difficult to
construct a sequence of step functions $U_2(x_2), U_3(x_2), \cdots,
U_N(x_2)$, so that $U_2$ is $L^2$-close to $U_1=U$, $U_N$ is
$L^2$-close to $V$, and every profile $U_k$ is obtained from $U_{k-1}$
by one of two operations, described in Theorems 3 and 4. Using these
theorems and the notes above, we construct a piecewise-smooth force
$f(x,t)$, such that $\int_0^T\parallel
f(\cdot,t)\parallel_{L^2}dt<\epsilon$, and $f$ transfers $U$ into $V$
during the time interval $[0,T]$.

{\bf 4. } Theorems 3 and 4 are proven by the variational method. 

Let $\cal D(M)=\cal D$ be the group of the volume-preserving
diffeomorphisms of the flow domain $M$. These diffeomorphisms may be
identified with fluid configurations: every configuration is obtained
from some fixed one by a permutation of fluid particles, which is
assumed to be a smooth, volume preserving diffeomorphism. The flow is
a family $g_t$ of elements of $\cal D$, depending on time $t, \ 0\le
t\le T$. The Lagrangian velocity of the flow is a vector-function
$V(x,t)=\fr{\pa}{\pa t}g_t(x)=\dot g_t(x)$, while the Eulerian
velocity is the vector field $v(x,t)=\dot g_t(g_t^{-1}(x))$. The action
of the flow is defined as $J\{g_t\}_0^T=\int_0^T\fr{1}{2}\parallel\dot
g_t\parallel_{L^2}^2 dt$, and the length
$L\{g_t\}_0^T=\int_0^T\parallel\dot g_t\parallel_{L^2} dt$.

The solution $u(x,t)$ of the homogeneous Euler equations (1), (2) is
an Eulerian velocity field of a geodesic trajectory $g_t$ on the group
$\cal D$: $u(x,t)=\dot g_t(g_t^{-1}(x))$, such that $\delta
J\{g_t\}_0^T=0$, provided $g_0, \ g_T$ are fixed. This implies that
also $\delta L\{g_t\}_0^T=0$. This is the classical Hamiltonian
principle (see[A3]). The evident idea is to try to construct solutions
of the Euler equations by fixing $g_0, \ g_T\in\cal D$, and looking
for the shortest trajectory, connecting these fluid configurations. If
the minimum is attainable, then we have constructed some nontrivial
solution of the Euler equations. 

But this idea does not work well. If $g_0$ and $g_T$ are $C^2$-close,
the minimum is assumed at some smooth trajectory.  But if $g_0$ and
$g_T$ are far away from each other, which is the only interesting
case, then it is possible that the minimum is no more attainable
(see[S], [A-K]). In the 3-d case there are examples of $g_0, g_T$,
such that for every smooth flow $g_t$, connecting $g_0$ and $g_T$,
there exists another smooth flow $g_t'$, connecting the same fluid
configurations, such that $J\{g_t'\}_0^T<J\{g_t\}_0^T$; so, the
minimum is unattainable. 

The existence of a minimal geodesic connecting two configurations of a
2-d fluid is neither proven nor disproven, while some physical
considerations show that sometimes the minimal smooth flow does not
exist . 

If there is no smooth solution of the variational problem, we may look
for a generalized solution, which is no longer a smooth flow, but
belongs to a wider class of object. The appropriate notion of a
generalized flow was introduced by Y. Brenier [B]. Generalized flow is
a probability measure $\mu$ in the space $X=C(0,T; M)$ of all
continuous
trajectories in the flow domain $M$, satisfying the following two
conditions:

\medskip

1. For every $t\in[0,T]$ and every Borel set $A\subset M$,

$$\mu\{x(\cdot)|x(t)\in A\}=\mes A;$$

2. 

$$J\{\mu\}=\int_X\int_0^T|\dot x(t)|^2 dt\mu\{dx\}<\infty$$

The meaning of the first condition is that the generalized flow is
incompressible; the second condition expresses the finiteness of
the mean action (and that $\mu$-almost all trajectories belong to
$H^1$). 

Every smooth flow may be regarded as a generalized one; but there is a
lot of truly generalized flows.

The generalized variational problem may be posed as follows: given a
diffeomorphism $g\in{\cal D}$; consider all generalized flows which,
in addition to the above conditions, satisfy the following:

\medskip

3. For $\mu$-almost all trajectories $x(t)$, \ $x(T)=g(x(0))$.

\medskip

We are looking for a generalized flow, which satisfies all three
conditions and minimizes the functional $J\{\mu\}_0^T=\int_X
\int_0^T\fr{1}{2}|\dot x(t)|^2 dt\mu\{dx\}$. 

Y.Brenier has proved, using the simple ideas of weak compactness of a
family of measures and semicontinuity of the action functional in $X$,
that this problem has a solution for every $g\in{\cal D}$. Simple
examples show that this solution may be very far from any smooth, or
even measurable, flow.

But in the 2-d case the situation is much better, because there is an
additional structure. To see it, consider a smooth incompressible flow
$g_t, \ 0\le t\le T, \ g_0=\id, \ g_T=g$. Let $Q=M\times[0,T]$ be a
cylinder in the $(x,t)$-space. Every trajectory
$\lambda_x=\{(g_t(x),t)\}, \
x\in M$, is a smooth curve in $Q$, connecting the points $(x,0)$ and
$(g(x),T)$. For different points $x, x'$, the curves $\lambda_x,
\lambda_{x'}$ do
not intersect. So, the lines $\lambda_x$ form a {\it braid},
containing
continuum of threads. Such a braid is called a smooth braid. 

Now let us define a generalized braid. Let us fix a volume-preserving
diffeomorphism $g\in{\cal D}$ and a piecewise-smooth incompressible
flow $G_t$,
such that $G_0=\id, \ G_T=g$. The bundle of lines $(G_t(x),t), \ x\in
M$, is called a reference braid and denoted by ${\bf B}_0$.

\medskip {\bf Definition 1. }{\sl A generalized flow $\mu$ is called a
generalized braid, and denoted by ${\bf B}$, if it satisfies
conditions 1, 2, 3 above, and the following condition

\medskip

4. For any $N$, let us pick trajectories $x^1(t), \cdots, x^N(t)$ by
random and independently, i.e. with a probability distribution
$\mu\otimes\cdots\otimes\mu$ ($N$ times). Then for almost all such
$N$-tuples of trajectories, the lines $(x^i(t),t)$ for different $i$ 
do not intersect, and the finite braid, formed by the curves
$(x^1(t),t), \cdots, (x^N(t),t)$, is isotopic to the braid, formed by
the curves $(G_t(x^1(0)),t), \cdots, (G_t(x^N(0)),t)$ (these braids
have the same endpoints, so it is possible to define their isotopy).}

\medskip

The braid ${\bf B}$ is called a braid weakly isotopic to a
piecewise-smooth braid ${\bf B}_0$. 

We are discussing the following variational problem: given a map $g$
and a reference braid ${\bf B}_0$; find a generalized braid
${\bf B}$, isotopic to ${\bf B}_0$, which minimizes the functional
$J\{{\bf B}\}$. 

\medskip

{\bf Theorem 5. }{\sl The variational problem has a solution for every
data \ $g, {\bf B}_0$.}

\medskip

To prove this theorem, we consider a sequence $\b_i$ of braids, such
that $J\{\b_i\}_0^T\searrow J_0$, where $J_0=\inf J\{\b\}$ for all
braids $\b$, isotopic to a given braid $\b_0$. This sequence is weakly
compact; its subsequence converges to a generalized flow $\mu$,
such that $J\{\mu\}=J_0$, exactly as for the generalized flows. But
the generslized flow $\mu$ is, in fact, a braid isotopic to $\b_0$,
which we shall denote by $\bar\b$.
This is implied by the fact that the isotopy class of a finite
subbraid with fixed endpoints is ''weakly continuous'', and therefore
we can pass to the limit and conclude that the weak limit of the
braids $\b_i$, regarded as generalized flows, is  a braid,
isotopic to $\b_0$. Let us call $\bar\b$ a minimal braid.

Generally, braids are as nonregular locally as generalized flows. In
particular, they have, in general, no definite velocity field: for
almost every $(x,t)$ there are different trajectories passing
through this point with different velocities. But  the {\it minimal}
braid is much more regular. Recall that a measurable flow is defined
as a family $h_{s, t}$ of measurable maps of $M$ into itself,
preserving the Lebesgue measure, and such that $h_{s, t}\circ h_{r,
s}=h_{r, t}$ for all $r, s, t \in [0,T]$.

\medskip

{\bf Theorem 6. }{\sl Let $\bar\b$ be a minimal braid, isotopic to the
reference braid $\b_0$. Then there exists
a measurable flow $h_{s, t}$ in $M$, such that for $\mu$-almost all
trajectories $x(t)$, \ $x(t)=h_{s, t}x(s)$. Moreover, there exists a
vector field $u(x,t)\in L^2$, divergence free and tangent to $\pa M$,
such that for almost all trajectories $x(t)$, \ $\dot x(t)=u(x(t),t)$
for almost all $t$.}

\medskip

The next fact about the minimal braids is the following

\medskip

{\bf Theorem 7. }{\sl The velocity field $u(x,t)$, corresponding to
the minimal braid $\bar\b$, is a weak solution of the Euler equations.
This means that for every vector field $v(x,t)\in C_0^{\infty}$, such
that $\nabla\cdot v=0$, and for every scalar function $\phi(x,t)$,

\be
&&\intl_Q\big[(u,\fr{\pa v}{\pa t})+(u\otimes u, \nabla v)\big] dx
dt=0,\\
&&\intl_Q (u,\nabla\phi) dx dt=0.
\ee
}

\medskip

The last fact which we need is the following approximation theorem.

\medskip

{\bf Theorem 8. }{\sl Suppose that $\bar\b$ is a minimal braid, and
$u(x,t)$ is its velocity field. Then for every $\epsilon>0$ there
exists a smooth incompressible flow with velocity field $w(x,t)$,
which is a solution of the nonhomogeneous Euler equations with the
force $f(x,t)$, such that $\parallel
w(x,t)-u(x,t)\parallel_{L^2}<\epsilon$ for all $t$, and
$\int_0^T\parallel f(\cdot,t)\parallel_{l^2}dt<\epsilon$.}

\medskip

{\bf 5. }The brais are used to construct the flows, described in
Theorems 3 and 4. Consider a flow with a piecewise-constant
profile $U(x_2)$; then the cylinder $Q=M\times[0,T]$ may be divided
into slices $Q_k$, so that $U|_{Q_k}=U^{(k)}$. In every such slice
the trajectories are parallel lines with the same slope. Thus, they
form a simple, piecewise-smooth braid.

Now let us describe the braids corresponding to the flows described in
Theorem 3. Let us divide the domains $M_k$, $M_{k+1}$, the bases of
the cylinders $Q_k$, $Q_{k+1}$, into small subdomains $M_{k,j}$,
$M_{k+1, l}$. Let us pick one point $x_{k,j}$ in every domain
$M_{k,j}$, and one point $x_{k+1, l}$ in every domain $M_{k+1, l}$.
Let $\lambda_{k,j}, \ \lambda_{k+1, l}$ be the trajectories, passing
through the points $(x_{k,j},0)$ and $(x_{k+1, l},0)$. Their endpoints
in $M\times\{T\}$ are denoted by $(y_{k,j},T)$ and $(y_{k+1, l},T)$.
The
trajectories, passing through $M_{k,j}\times \{0\}$, and through
$M_{k+1, l}\times\{0\}$, form subbraids $\b_{k,j}$ and $\b_{k+1, l}$.

Now let us define a new braid $\b'_0$. First let us define its threads
$\lambda'_{k,j}$ and $\lambda'_{k+1, l}$, passing through the points
$(x_{k,j},0)$ and $(x_{k+1, l},0)$. They are straight lines, passing
through
the points $(y'_{k,j},T)$ and $(y'_{k+1, l},T)$, obtained from the
points
$(y_{k,j},T)$ and $(y_{k+1, l},T)$ by the shift in the $x_2$-direction
by, resp., $(x_2^{(k+1)}-x_2^{(k)})$ and $(x_2^{(k-1)}-x_2^{(k)})$.
These lines do not, generally, intersect.

Now let us define a piecewise-smooth braid $\b'_0$. It coincides with
$\b$ outside $Q_k\bigcup Q_{k+1}$. In the last domain the braid $b'_0$
consists of smooth incompressible subbraids $\b'_{k,j}$ and $\b'_{k+1,
l}$ with the bases $M_{k,j}\times\{0\}$ and $M_{k+1, l}\times\{0\}$.
Each subbraid contains one line $\lambda'_{k,j}$ and $\lambda'_{k+1,
l}$. The interfaces between these subbraids are piecewise-smooth. It
is easy to construct such subbraids, while they are not unique.

Now let us use $\b'_0$ as a reference braid, and construct a
minimal braid $\b'$, isotopic to $\b_0'$. Using Theorems 5--8, we
construct a smooth flow $w(x,t)$, supported by a smooth force
$f(x,t)$, such that $\int_0^T\parallel f(\cdot,t)\parallel_{L^2}dt$ is
arbitrarily small, provided $T$ is big enough. This flow is
$L^2$-close to $U$ at $t=0$ and to $V$ at $t=T$; after a small
modification of $w(x,t)$, requiring an $L^2$-small correcting force,
we obtain a flow described in Theorem 3.

The proof of Theorem 4 is similar; but in this case the curves
$\lambda'_{k,j}$ and $\lambda'_{k+1, l}$ \  are going from points
$(x_{k,j},0) \ ( (x_{k+1, l},0) )$ to the points $(y_{k,j},0) \ (
(y_{k+1, l},0) )$ , and some of the curves $\lambda'_{k,j}$ and
$\lambda'_{k+1, l}$ are linked. 

{\bf 6. } Theorems 3 and 4 are true also for circular flows in a disk,
with the angular momentum staying in place of momentum in Theorem 4.
But for generic 2-d domains the situation is not so clear. We don't
know, whether there is an integral of motion, similar to the angular
momentum, in any domain different from the disk. If such integral does
not exist, which is the most likely, then the natural conjecture is
that for any two flows with equal energies the conclusion of Theorem 4
is true. But this behavior is paradoxical: just imagine a nearly
circular flow in a nearly circular domain (e.g. ellipse), which after
some long time changes the sign of the angular velocity. This question
requires more thinking.

\medskip

{\bf Acknowledges. } This work was done during my stay at the
Department of Mathematics of Princeton University in the spring
semester 1999, with the support of the American Institute of
Mathematics. I am very thankful to these institutions for creative
atmosphere and generous support.

I am thankful to the organizers of the Journees des Equations
Differentielles for invitation to this excellent conference.

\bigskip
\bigskip

\centerline{LIST OF REFERENCES}
\medskip

[A1] V. Arnold, {\it Sur la G\' eom\' etrie diff\' e rentielle des
groupes de Lie de dimension infinie et ses applications \` a
l'hydrodynamique des fluides parfaits}, Ann. Inst. Fourier {\rm 16}
(1966), 316-361.

[A2] V. Arnold, {\it On the a priori estimate in the theory of
hydrodynamical stability}, Amer. Math. Soc. Transl. {\bf 19} (1969),
267-269.

[A3] V. Arnold, {\it Mathematical methods of classical mechanics},
Springer-Verlag, New York, 1989.

[A-K] V. Arnold, B. Khesin, {\it Topological methods in
hydrodynamics}, Applied Mathematical Sciences, v. 125, 
Springer-verlag, 1998.

[B] Y. Brenier, {\it The least action principle and the related
concept of generalized flows for incompressible perfect fluids}, J.
Amer. Math. Soc. {\bf 2} (1989), no. 2, 225-255.

[M-P] C. Marchioro, M. Pulvirenti, {\it Mathematical theory of
incompressible nonviscous fluids}, Applied Mathematical Sciences, v.
96, Springer-Verlag, 1994.

[S] A. Shnirelman, {\it The geometry of the group of diffeomorphisms
and the dynamics of an ideal incompressible fluid}, Math. USSR Sbornik
{\bf 56} (1987), no. 1, 79-105.

\end{document}